\documentclass[12pt]{amsart}
\let\oldlabel=\label
\def\prellabel{\marginparsep=1em\marginparwidth=44pt
     \def\label##1{\oldlabel{##1}\ifmmode\else\ifinner\else
          \marginpar{{\footnotesize\ \\ \tt
                     ##1}}\fi\fi}}
\def\NN{{\NZQ N }}

\def\ZZ{{\NZQ Z}}
\def\RR{{\NZQ R}}

\def\opn#1#2{\def#1{\operatorname{#2}}} 
\opn\chara{char}
\opn\gr{gr}
\opn\rank{rank}
\let\sect=\cap

\let\iso=\cong

\let\Sect=\bigcap
\let\Dirsum=\bigoplus

\let\:=\colon
\newtheorem{theorem}{Theorem}[section]
\newtheorem{lemma}[theorem]{Lemma}
\newtheorem{corollary}[theorem]{Corollary}
\newtheorem{proposition}[theorem]{Proposition}
\theoremstyle{definition}
\newtheorem{remark}[theorem]{Remark}

\let\epsilon=\varepsilon
\let\phi=\varphi
\let\kappa=\varkappa
\opn\ini{in}
\opn\KRS{KRS}

\def\pp{{\mathfrak p}}
\def\qq{{\mathfrak q}}
\opn\krs{krs}
\opn\diag{diag}
\opn\DD{{\mathcal D}}
\opn\SS{{\mathcal S}}
\opn\MM{{\mathcal M}}
\opn\GL{GL}

\def\XX{{\mathcal X}}
\def\Rees{{\mathcal R}_t}
\let\Bbb=\mathbb

\def\RR{{\Bbb R}}

\def\ZZ{{\Bbb Z}}
\def\NN{{\Bbb N}}
\def\ini{\operatorname{in}}

\opn\height{height}
\opn\length{length}
\opn\cl{cl}
\opn\Cl{Cl}
\opn\Grass{Grass}

\begin{document}

\title{ Algebras of minors }
\author{Winfried Bruns \and Aldo Conca}
\address{Universit\"at Osnabr\"uck,
FB Mathematik/Informatik, 49069 Osnabr\"uck, Germany}
\email{Winfried.Bruns@mathematik.uni-osnabrueck.de}
\address{Dipartimento di Matematica e Fisica
Universit\'a di Sassari, Via Vienna 2, I-07100 Sassari, Italy }
\email{conca@ssmain.uniss.it}

\dedicatory{To J\"urgen Herzog on the occasion of his 60th
birthday}

\maketitle

\section*{Introduction}
Let $I$ be an ideal of a Noetherian ring $R$. The Rees algebra of
$I$ is the graded $R$-algebra $\Dirsum_{k=0}^{\infty} I^k$. The
study of the properties of Rees algebras (and of the other blow-up
algebras) has attracted the attention of many researchers in the
last three  decades. For a detailed account on  the theories that
have been developed and  on the results that have been proved  the
reader should consult the monograph of Vasconcelos \cite{V}.  In
this paper we treat a special and interesting  case: the Rees
algebras of determinantal ideals and their special fibers.

Let $X=(x_{ij})$ be a generic matrix of size $m\times n$ over a
field $K$ and let $S$  be the polynomial ring $K[x_{ij}]$. Let
$I_t$ be the ideal of $S$ generated by the minors of size $t$ of
$X$. Finally, let $\Rees$ be the Rees algebra  of $I_t$, and let
$A_t$ be the subalgebra of $S$ generated by the $t$-minors of $X$.
In the case of maximal minors, i.e.\ $t=\min(m,n)$, the algebra
$A_t$ is nothing but the homogeneous coordinate ring of the
Grassmann variety which is known to be a Gorenstein factorial
domain; for instance, see Bruns and Vetter \cite{BV}. Furthermore,
for maximal minors, Eisenbud and Huneke proved in \cite{EH} that
$\Rees$ is normal and Cohen-Macaulay. Their approach is based on the notion
of algebra with straightening law. Later Bruns proved that $\Rees$
and $A_t$ are Cohen-Macaulay and normal for any $t$ if $\chara
K=0$ using invariant theory methods; see \cite{B}.  In \cite{BC}
we have shown, by using Sagbi basis deformation and the
Knuth-Robinson-Schensted correspondence, that $\Rees$ and $A_t$
are  Cohen-Macaulay and normal domains for any $t$ and for any
non-exceptional characteristic. Recall that (in our terminology)
$K$ has non-exceptional characteristic if either $\chara K=0$ or
$\chara K>\min(t,m-t,n-t)$. For exceptional characteristic $\Rees$
and $A_t$ can be (and perhaps are always) very far from being
Cohen-Macaulay; see the example in \cite{B}.

For maximal minors the divisor class group and the canonical class
of $\Rees$ have been determined. They are the ``expected" ones,
that is, those of the Rees algebra of a prime ideal   with primary
powers in a regular local ring: $\Cl(\Rees)$ is free of rank $1$
with generator $\cl(I_t\Rees)$ and the canonical class is
$(2-\height I_t)\cl(I_t\Rees)$; see Herzog and Vasconcelos
\cite{HV} or Bruns, Simis and Trung \cite{BST}. On the other hand,
for non-maximal minors, $\Cl(\Rees)$ is free of rank $t$ \cite{B}.

Our goal is to determine the canonical class of $\Rees$, the
divisor class group of $A_t$ and the canonical class of $A_t$ in
all the remaining cases, that is, $t<\min(m,n)$.  The case $t=1$
is anyway easy, and also the study of $A_t$ in the case $m=n=t+1$
is trivial since that ring is a polynomial ring.

Our main tools are the $\gamma$-functions that allow us to
describe all relevant ideals in $\Rees$ and $A_t$ and the Sagbi
deformation by which we will lift the canonical module from the
initial algebras of $\Rees$ and $A_t$ to the algebras $\Rees$ and
$A_t$ themselves. In the following we consider $A_t$ and $\Rees$
as subalgebras of the polynomial ring $S[T]$, $A_t$ being
generated by the elements of $M_t T$, where $M_t$ is the set of
$t$-minors, and $\Rees$ being generated by $M_t T$ and the entries
of $X$.

It turns out that the canonical class of $\Rees$ is given by
\begin{align*}
\cl(\omega(\Rees))&=\sum_{i=1}^t \Bigl(2-(m-i+1)(n-i+1)+t-i\Bigr)
\cl(P_i)\\
&=\cl(I_t\Rees)+\sum_{i=1}^t (1-\height I_i)\cl(P_i)
\end{align*}
Here $P_i$ is the prime ideal of $\Rees$ determined by $\gamma_i$;
it has the property $P_i\sect S=I_i$, and it is well-known that
$\height I_i=(m-i+1)(n-i+1)$. Moreover,
$\cl(I_t\Rees)=\sum_{i=1}^t (t-i+1)\cl(P_i)$. The classes
$\cl(P_i)$ generate $\Cl(\Rees)$ freely.

The divisor class group of $A_t$ is also free; its rank is $1$
with generator the class $\cl(\qq)$ of the ideal $\qq=(f)S[T]\sect
A_t$ where $f$ is a minor of size $t+1$ of $X$. The canonical
class of $A_t$ is given by
$$
(mn-mt-nt)\cl(\qq).
$$
It follows that, apart from the trivial or known cases listed
above, $A_t$ is Gorenstein if and only if $mn=t(m+n)$.

There are many open questions concerning the algebras $\Rees$ and
$A_t$. Just to mention few of them:\smallskip

\noindent (a) In general we do not know (the degrees of) the
defining equations of $\Rees$ and $A_t$. It is known that $\Rees$
and $A_t$ are defined by quadrics in the maximal minors case or if
$t=1$ or if $m=n=t+1$  and that in general higher degree relations
are needed.
Since $A_t$ has a rational singularity, it has a negative $a$-invariant.
This implies that the defining equations of $A_t$ have degree $\leq \dim
A_t$ and $\dim A_t=mn$  for non-maximal minors. But this bound is not at
all sharp.
There are indications that quadrics and cubics are
enough for $A_t$ and $\Rees$ if $t=2$.
\smallskip

\noindent (b)  The Hilbert function of $A_t$ can be computed by
the hook length formula (for example, see Stanley \cite[7.21.6]{St}). We
have in fact
computed these Hilbert functions for a substantial number of
cases, but we do not know an explicit formula, not even, say, for
the multiplicity of $A_t$.
\smallskip

\noindent (c)   For maximal minors  the algebra  $A_t$ has an
isolated singularity since the Grassmann variety is a homogeneous
variety and, hence, smooth; for instance see \cite{BV}. But in
general we do not know (the dimension of) the singular locus of
$A_t$. \smallskip

In principle, the results of this paper can be extended to the
other generic determinantal ideals (minors of a symmetric matrix,
pfaffians of an alternating matrix, etc.).  In the last section we
do this in detail for a generic Hankel matrix.

Most of the this work was done during the visit of the first
author to the University of Genova  in the frame of a project
supported by the Italian Indam-Gnsaga.

\section{Notation and generalities}

We keep the notation of the introduction.  We will always assume
that the characteristic of the field $K$ is not exceptional, which
means that it is either $0$ or larger than $\min(t,m-t,n-t)$. We
will study the algebras $\Rees$ and $A_t$. As mentioned above,
both algebras can be embedded into the polynomial ring $S[T]$
where $T$ is an auxiliary indeterminate. Denote by $M_t$ the set
of the $t$-minors of $X$. Then $A_t$ is (isomorphic to) the
$K$-subalgebra of $S[T]$ generated by the elements of the set
$M_tT$, and $\Rees$ is the subalgebra generated by the elements of
$X$  and of $M_tT$. Clearly, $A_t\subset \Rees$ and both algebras
are $\NN$-graded in a natural way. Moreover $\Rees/(X) \iso A_t$
where $(X)$ is the ideal of $\Rees$ generated by the $x_{ij}$.
Denote by  $V_t$ the $K$-subalgebra of $S[T]$ generated by
monomials that have degree $t$ in the variables $x_{ij}$ and
degree $1$ in $T$. Note that $A_t\subset V_t$ and that $V_t$ is a
normal semigroup ring isomorphic to the $t$-th Veronese subring of
the polynomial ring $S$. The rings $\Rees$ and $A_t$ are known to
be normal Cohen-Macaulay domains. The  goal of the paper is to
determine their divisor class group and canonical class and to
discuss the Gorenstein property of these rings. As we have
mentioned in the introduction these invariants are already known
is some special cases. Therefore we will concentrate our attention
to the remaining cases, that is, we make the following
\smallskip

\noindent{\bf Assumptions}.  (a) When studying the Rees algebra
$\Rees$, we will assume  that  $t<\min(m,n)$, \\[3pt]
\noindent (b) when studying $A_t$, we will assume that
$1<t<\min(m,n)$ and that $m\neq n$ if $t=\min(m,n)-1$.
\smallskip

Our approach to the study of the algebras under investigation
makes use of Sagbi bases deformations and of the straightening law
for generic minors. For generalities  on the former we refer the
reader to Conca, Herzog and Valla \cite{CHV}.  As for the latter,
the reader can consult Bruns and Vetter \cite{BV}, Doubilet, Rota and Stein
\cite{DRS}, or De Concini, Eisenbud and Procesi \cite{DEP}.

Here we just recall its main features.  The minor of $X$ with row
indices $a_1,\dots,a_t$ and column  indices $b_1,\dots,b_t$ is
denoted by
$$
[a_1,\dots,a_t | b_1,\dots,b_t].
$$
We will assume that the row and the column indices are given in
ascending order $a_i<a_{i+1}$ and  $b_i<b_{i+1}$. The set of the
minors of $X$ is partially ordered in a natural way:
$$
[a_1,\dots,a_t | b_1,\dots,b_t]\leq [c_1,\dots,c_s | d_1,\dots,d_s]
$$
if and only if $t\geq s$ and $a_i\leq c_i$, $b_i\leq d_i$ for all
$i=1,\dots,s$.

Let  $\delta_1,\delta_2,\dots,\delta_p$ be minors of size
$t_1,t_2,\dots,t_p$ respectively.  Consider the product
$\Delta=\delta_1\cdots \delta_p$. The shape of  $\Delta$ is  by
definition the sequence of numbers  $t_1,t_2,\dots,t_p$ and
$\Delta$  is said to be a standard monomial if $\delta_1\leq \dots
\leq \delta_p$.  The straightening law asserts the following:

\begin{theorem}\label{stra}
The standard monomials form a $K$-vector space basis of $S$.
\end{theorem}

We recall now the definition of the functions $\gamma_t$. Given a
sequence   of numbers  $s=s_1,\dots,s_p$ and a number $t$ we
define
$$
\gamma_t(s)=\sum_{i=1}^p \max(s_i+1-t,0).
$$

Then we extend the definition to products of minors by setting:
$$
\gamma_t(\Delta)=\gamma_t(s)
$$
where $s$ is the shape of $\Delta$. Finally,  for ordinary
monomials $M$ of $S$ we define:
$$
\rho_t(M)=\sup\{ \gamma_t(\Delta) :   \Delta \text{  is a product
of minors of $X$ with } \ini(\Delta)=M\}
$$
where $\ini(f)$ denotes the initial term of a polynomial $f$ with
respect to a diagonal term order, i.e.\ a term order such that
$\ini(\delta)=x_{a_1b_1}\cdots x_{a_tb_t}$ for every minor
$\delta=[a_1,\dots,a_t| b_1,\dots,b_t]$ of $X$.  In our previous
papers \cite{BC,BC1}, the function $\rho_t$  was also denoted by
$\gamma_t$; this could create some ambiguities in the present
paper, so that we have changed the notation.

The functions $\gamma_t$ were introduced in \cite{DEP} to describe
the symbolic powers of $I_t$ and the primary decomposition of the
powers of $I_t$. We have shown in \cite{BC} that the functions
$\rho_t$ describe instead the initial ideals of the symbolic and
ordinary powers of $I_t$. Let us recall the precise statements:

\begin{theorem}
\label{symbpow} \leavevmode\kern5pt
\begin{itemize}
\item[(a)] $I_t^{(k)}$ is  $K$-vector space generated by the standard
monomials  $\Delta$ with $\gamma_t(\Delta)\geq k$, and it contains
all products of minors $\Delta$ with $\gamma_t(\Delta)\geq k$.
\item[(b)] The initial ideal $\ini(I_t^{(k)})$  is the  $K$-vector space
generated by the ordinary monomials $M$ with $\rho_t(M)\geq k$.
\end{itemize}
\end{theorem}

\begin{theorem}
\label{prdec}\leavevmode\kern5pt
\begin{itemize}
\item[(a)]   $I_t^{k}=\Sect_{i=1}^t  I_i^{(k(t+1-i))}$ is a primary
decomposition of $I_t^k$. In particular,  $I_t^k$ is generated by
the standard monomials  $\Delta$ with $\gamma_i(\Delta)\geq
k(t+1-i)$ for all $i=1,\dots,t$, and it contains all products of
minors $\Delta$ satisfying these conditions. Moreover this primary
decomposition is irredundant if $t<\min(m,n)$ and $k\geq
(v-1)/(v-t)$ where $v=\min(m,n)$.
  \item[(b)] The initial ideal $\ini(I_t^{k})$  is the  $K$-vector space
generated by the ordinary monomials $M$ with $\rho_i(M)\geq
k(t+1-i)$ for all $i=1,\dots,t$.
\end{itemize}
\end{theorem}

We may define  the value of  function $\gamma_t$ for any
polynomial $f$ of $S$  as follows. Let $f=\sum_{i=1}^p \lambda_i
\Delta_i$ the unique representation of $f$ as a linear combination
of standard monomials $\Delta_i$ with coefficients  $\lambda_i\neq 0$. Then
we set
$$
\gamma_t(f)=\inf\{ \gamma_t(\Delta_i)  :   i=1,\dots,p\}
$$
and $\gamma_t(0)=+\infty$.  This definition is consistent with the
one above in the sense that if $\Delta$ is a product of minors of
shape $s$ which is non-standard, then $\gamma_t(\Delta)=
\gamma_t(s)$.

The function $\gamma_t$ is indeed a discrete valuation on $S$
(with values in $\NN$), that is, the following conditions are
satisfied for every $f$ and $g$ in $S$:
\begin{itemize}
\item[(a)] $\gamma_t(f+g)\geq \min(\gamma_t(f), \gamma_t(g))$, and
equality holds if $\gamma_t(f)\neq \gamma_t(g)$,
\item[(b)] $\gamma_t(fg)= \gamma_t(f)+\gamma_t(g)$.
\end{itemize}
Note that (a) follows immediately from the definition, and (b)
from the fact that the associated graded ring of a symbolic
filtration is a domain if the base ring is regular. We
may further extend $\gamma_t$ to the field of fractions $Q(S)$ of
$S$ by setting
$$
\gamma_t(f/g)=\gamma_t(f)-\gamma_t(g)
$$
so that $S$ is a subring of the valuation ring associated with each
$\gamma_t$.

The next step is to extend the valuation $\gamma_t$ to the
polynomial ring $S[T]$ and to its field of fractions. We want to
do this in a way such that the subalgebras $\Rees$ and $A_t$ of
$S[T]$ will then be described in terms of these functions.  So,
from now on, let us fix a number $t$ with $1\leq t\leq \min(m,n)$.
For every polynomial
$F=\sum_{j=0}^p f_jT^j\neq 0$ of $S[T]$ we set
$$
\gamma_i(F)=\inf\{\gamma_i(f_j)-j(t+1-i):j=0,\dots,p\};
$$
in particular
$$
\gamma_i(T)=-(t+1-i).
$$

Then we have:

\begin{proposition} \label{gaQ}
The function $\gamma_i$ defines a discrete valuation on the field
of fractions $Q(S[T])$ of $S[T]$.
\end{proposition}

\begin{proof}
This is a general fact, see Bourbaki \cite[Ch.\ VI, \S10, no.\ 1,
lemme 1]{Bo}.
\end{proof}

Note that $F=\sum f_jT^j\in S[T]$ has $\gamma_i(F)\geq 0$ if and
only if $f_j\in I_i^{(j(t+1-i))}$ for every $j$.  It follows from
\ref{prdec} that the Rees algebra $\Rees$ of $I_t$ has the
following description:

\begin{lemma} \label{gaRees}
$$
\Rees= \{ F \in S[T] :  \gamma_i(F)\geq 0\text{ for }i=1,\dots,
t\}.
$$
\end{lemma}

Similarly $A_t= \{ F \in V_t : \gamma_i(F)\geq 0\text{ for }
i=2,\dots, t\}$.  But this
description is redundant:

\begin{lemma}  \label{gaA}
$$
A_t= \{ F \in V_t : \gamma_2(F)\geq 0  \}.
$$
\end{lemma}

\begin{proof} We have to show that if $\Delta$ is a standard monomial
with $\deg \Delta=kt$ and $\gamma_2(\Delta)\geq k(t-1)$, then
$\gamma_i(\Delta)\geq k(t+1-i)$ for all $j=3,\dots,t$. For $i=1,
\dots, m$ let $a_i$ denote the number of the factors of $\Delta$
which are minors of size $i$.  By assumption we have that
$$
\text{(a)} \ \   \sum_{i=1}^m ia_i=kt  \qquad\text{and} \qquad
\text{(b)} \ \ \sum_{i=2}^m (i-1)a_i\geq k(t-1).
$$
In view of (a), condition (b) can be rewritten as
$$
\text{(c)}\ \   \sum_{i=1}^m a_i\leq k.
$$
We have to show that
$$
\text{(d)}\ \ \sum_{i=j}^m a_i (i+1-j)\geq k(t+1-j).
$$
Note that
\begin{align*}
\sum_{i=j}^m a_i (i+1-j)&=\sum_{i=1}^m a_i
(i+1-j)-\sum_{i=1}^{j-1} a_i (i+1-j)\\&=kt+\sum_{i=1}^m a_i
(1-j)+\sum_{i=1}^{j-1} a_i (j-i-1).
\end{align*}
Therefore (d) is equivalent to
$$
\biggl(k-\sum_{i=1}^m a_i\biggr) (j-1)+\sum_{i=1}^{j-1} a_i
(j-i-1)\geq 0
$$
which is true since the left hand side is the sum of non-negative
terms.
\end{proof}

Similarly the  initial algebras of $\Rees$ and $A_t$ have a
description in terms of the functions $\rho_i$. To this end we
extend the definition of the function $\rho_i$ to  monomials of
$S[T]$ by setting
$$\rho_i(MT^k)=\rho_i(M)-k(t+1-i)$$
where $M$ is a monomial of $S$. Furthermore, for a polynomial
$F=\sum_{j=1}^p  \lambda_j N_j$  of $S[T]$ where the $N_i$ are
monomials and the $\lambda_j$ are non-zero elements of $K$, we set
$$
\rho_i(F)=\inf \{ \rho_i(N_j) : j=1,\dots,p\}.
$$
Now Theorem \ref{prdec} implies:

\begin{lemma} \label{gaIn}
\leavevmode\kern5pt
\begin{itemize}
\item[(1)]
The initial algebra $\ini(\Rees)$ of $\Rees$ is given by
$$\ini(\Rees)=\{ F \in S[T] :  \rho_i(F)\geq 0  \mbox{ for  } i=1,\dots,t\}$$
The set of the monomials $N$ of $S[T]$ such that $\rho_i(N)\geq 0$
for $i=1,\dots,t$ form  a $K$-vector space basis of $\ini(\Rees)$.

\item[(2)] The initial algebra $\ini(A_t)$ of $A_t$ is given by:
$$\ini(A_t)=\{ F \in V_t :  \rho_2(F)\geq 0 \}$$
The set of the monomials $N$ of $S[T]$ such that $\rho_1(N)=0$ and
$\rho_2(N)\geq 0$ form  a $K$-vector space basis of $\ini(A_t)$.
\end{itemize}
\end{lemma}

The major difference between the functions $\gamma_i$ and $\rho_i$ is that
the latter  is not a valuation. Nevertheless, it will turn out that
$\rho_i$ is an ``intersection'' of valuations, see Section \ref{canoini}.

\section{The divisor class groups of $\Rees$ and $A_t$}

The divisor class group of $\Rees$ has already been determined in
\cite[Cor.\ 2.4]{B}. In this section we will present a (slightly)
different approach to it, which will also be used to determine the
divisor class group of $A_t$. We will need the following general facts:

\begin{lemma}
\label{referee}
Let $A$ be a normal  domain. Suppose $v_1,\dots,v_t$ are discrete valuations on
$A$.  Let $B=\{x\in A:  v_i(x)\ge 0 \mbox{ for } i=1,\dots,t\}$
and $P_i=\{x\in B: v_i(x)\ge 1\}$. Then:
\begin{itemize}
\item[(1)]  Assume that there exist elements  $y_1,\dots,y_t\in A$
with $v_i(y_i)<0$ and
$v_j(y_i)\geq 0$ for $j\neq i$.
  Then $P_i$ is a height one ideal of the Krull domain
$B$ whose symbolic power $P_i^{(j)}$ is $\{x\in B:v_i(x)\ge j\}$ for
every $j\ge 1$.
\item[(2)] Assume that $\height(P_iA) >1$ for    $i=1,\dots,t$ . Then
$A = \bigcap B_Q,$ where the intersection is taken over all height
1 prime ideals $Q$ of $B$ with $Q\not\in \{P_1,\dots,P_t\}$.
\end{itemize}
\end{lemma}

\begin{proof}
(1)  Let $V_i$ be the valuation ring of $v_i$. Then evidently
$$
B=V_1\cap \dots \cap V_t \cap \bigcap_{\height\ Q=1}A_Q,
$$
and $B$ is obviously a Krull domain. All essential valuation
overrings of $B$ occur in the above intersection. Since, by
assumption, $y_i$ is in $A$ and  in all the $V_j$ with $j\neq i$,
but not in $B$, we have that  $V_i$ is not redundant in the above
description of $B$. Then $V_i$ is an  essential valuation overring
of $B$. The rest is clear by standard results about essential
valuations of Krull domains.

(2)   Let $f$ be an arbitrary element of $A$. Let $Q$  be a height
1 prime ideal of $B$ with $Q\not\in \{P_1,\dots,P_t\}$.  We choose
an element $g \in P_1\cap \dots \cap P_t  \setminus Q$. Since
$v_i(g) \ge 1$ for all $i$,   there is a positive integer $r$ such
that $v_i(fg^r) = v_i(f) + rv_i(g) \ge 0$ for all $i$. It follows
that   $fg^r \in B$. Hence $f \in B_Q$. So we get $A \subseteq
B_Q$. To prove $A \supseteq \cap B_Q$, we first note that $A =
\bigcap A_{\mathfrak Q}$, where $\mathfrak Q$ is   taken over all
height 1 prime ideal $\mathfrak Q$ of $A$. So we only need to show
that if $Q = {\mathfrak Q} \cap B$, then $Q$ is a height 1 prime
ideal   of $B$ and $Q \not\in\{P_1,\dots,P_t\}$.  That  $Q
\not\in\{P_1,\dots,P_t\}$ follows from the assumption
$\height(P_iA) > 1$. To prove $\height(Q) = 1$, it is sufficient
to show that $B_Q = A_{\mathfrak Q}$. Let $f \in A\setminus
\mathfrak Q$. Then $f \not\in QB_Q$ because otherwise there is an
element $g \in B \setminus Q$ such that $fg \in Q \subset
\mathfrak Q$, contradicting the facts that $f \not\in \mathfrak
Q$, $g \not\in \mathfrak Q$. Thus, $f$ is an invertible element in
$B_Q$. Hence we can conclude that $A_{\mathfrak Q} \subseteq B_Q$
so that $A_{\mathfrak Q} = B_Q$.
\end{proof}

We are grateful to the referee for the above lemma, which has
simplified our original treatment. For $i=1,\dots,t$ we set
$P_i=\{ F \in \Rees : \gamma_i(F)\geq 1\}$.

\begin{lemma}
\label{symbpowR} The ideal $P_i$ is a height $1$ prime ideal of
$\Rees$ for $i=1,\dots, t$ and $P_i^{(j)}=\{ F \in \Rees :
\gamma_i(F)\geq j\}$ for every $j>0$.
\end{lemma}

\begin{proof} We apply \ref{referee}(1) with $A=S[T]$ and $v_i=\gamma_i$.  For
$i=1,\dots,t$ set
$$
y_i=g^af^{t-i}T^{m-i}
$$
where $g$ is an $(i-1)$-minor, $f$ is an $m$-minor and $a$ is a
large enough integer. The reader may check that these elements
have the right $\gamma$-values, i.e.\ $\gamma_j(y_i)\geq 0$  for
$j=1,\dots,t$ and $j\neq i$ and $\gamma_i(y_i)< 0$.
\end{proof}

We can now describe the polynomial ring as a subintersection of
the Rees algebra:

\begin{proposition} \label{subRees}
Set $R=\Rees$.  Then
$$
S[T]=\Sect R_Q
$$
where the intersection is extended over all the height $1$ prime
ideals $Q$ of $\Rees$ different from $P_1,\dots,P_t$.
\end{proposition}

\begin{proof} We apply \ref{referee}(2) with $A=S[T]$ and
$v_i=\gamma_i$.  To this end it
suffices to note that $I_iS[T] \subset P_iS[T]$ for all $i$.
\end{proof}

The ideal $I_t\Rees$ is a height $1$ ideal of $\Rees$. This
follows, for instance, from the fact that $\Rees$ has dimension
$\dim S+1$ and the associated graded ring $\Rees/I_t\Rees$ has dimension
$\dim S$. Furthermore $I_t\Rees$ is divisorial since it is
isomorphic to $\Dirsum_{k>0} I_t^kT^k$, which is a height $1$
prime ideal.

As a consequence of \ref{symbpowR} and \ref{prdec} (or \cite[Thm.\
2.3]{B}) we obtain

\begin{proposition}
\label{primdecRees}
One has
$$
I_t\Rees=\Sect_{i=1}^t  P_{i}^{(t-i+1)},
$$
and this is an irredundant primary decomposition of $I_t\Rees$. In
particular, $P_i$ is a height $1$ prime ideal of $\Rees$ for
$i=1,\dots,t$.
\end{proposition}

\begin{proof}
By virtue of \ref{symbpowR} the equality $I_t\Rees=\Sect_{i=1}^t
P_{i}^{(t-i+1)}$ is simply a re-interpretation of \ref{prdec}.
Since we assume that $t<\min(m,n)$, the primary decomposition of
$I_t^k$ given in \ref{prdec} is irredundant for $k\gg0$. It
follows that the primary decomposition of $I_t\Rees$ is also
irredundant.
\end{proof}

\begin{remark}
Since one can prove directly that $P_1$ is a height $1$ prime
ideal with primary powers, one can also use the standard
localization argument in order to show \ref{symbpowR} and
\ref{primdecRees} (see \cite{B}). However, this argument is not
available for the Hankel matrices discussed in Section
\ref{Hankel}.
\end{remark}

\begin{theorem}\label{divclRees}
The divisor class group $\Cl(\Rees)$ of $\Rees$
is free of rank $t$,
$$
\Cl(\Rees)\iso \ZZ^t
$$
with basis $\cl(P_1),\dots,\cl(P_t)$.
\end{theorem}

\begin{proof} The conclusion follows from \ref{primdecRees} and from a
general result of Simis and Trung \cite[Thm.\ 1.1]{ST}. Let us
mention that one can derive the result also directly from Fossum
\cite[Thm.\ 7.1]{F} and \ref{subRees}.
\end{proof}

We have similar constructions and results for $A_t$. First define
$$
\pp=\{F\in A_t : \gamma_2(F)\geq 1\}.
$$
It is clear that $\pp$ is a prime ideal of $A_t$.

\begin{lemma}\label{ht1At}
The ideal $\pp$ is prime of height $1$. Moreover, one has
$\pp^{(j)}=\{ F \in A_t : \gamma_2(F)\geq j\}$.
\end{lemma}

\begin{proof} We apply \ref{referee}(1) where $A=V_t$
and $v_1=\gamma_2$.  It suffices to take $f=x_{11}^tT\in V_t$ and
note that $\gamma_2(f)=-t+1<0$.
\end{proof}

\begin{proposition}\label{subintAt}
Set $A=A_t$. Then
$$
V_t=\Sect A_P
$$
where the intersection is extended over all the height $1$ prime
ideals $P$  of $A_t$  with $P\neq \pp$.
\end{proposition}

\begin{proof} We apply \ref{referee}(2) where $A=V_t$
and $v_1=\gamma_2$.  It suffices to prove that $\pp A$ has height
$>1$. To this end take $f_1$ and $f_2$ distinct $(t+1)$-minors
(they exist by our assumptions). Set  $g_j=f_j^tT^{t+1}$. Then
$\gamma_2(g_j)=1$ and hence $g_j\in \pp V_t$. Since the $f_i$ are
prime elements, $(f_1^t,f_2^t)S$ has height $2$. Since $V_t$ is a
direct summand of the polynomial ring it follows that
$(g_1,g_2)V_t$ has height $2$.
\end{proof}

Let $f$ be a $(t+1)$-minor of $X$. Set $g=f^tT^{t+1}$. By
construction, $g\in V_t$ and $\gamma_2(g)=1$ so that $g\in A_t$. Set
$$
\qq=(f)S[T]\sect A_t.
$$
In other words, $\qq=\{fa\in A_t  :  a \in S[T]\}$. Since $f$ is a
prime element in $S[T]$, the ideal $\qq$ is prime. Furthermore we
have $\pp\qq^t \subset (g) \subset \pp\sect \qq$. The second
inclusion is trivial. As for the first, note that any generator of
$\pp\qq^t$ can be written in the from $gb$,  and then just
evaluate $\gamma_2$ to show that $b$ is in $A_t$.

Note also that:

\begin{corollary}\label{fsectAt}
The ideal  $\qq$ is a height $1$ prime ideal of $A_t$. Furthermore
$\qq^{(j)}=(f^j)S[T]\sect A_t$.
\end{corollary}

\begin{proof}
It suffices to check that $\qq$ does not contain $\pp$. Let $f_1$
be a $(t+1)$-minor of $X$ different from $f$ (it  exist by our
assumptions). Then $f_1^tT^{t+1}$ is in $\pp$, but not
in $\qq$. The second statement is easy.
\end{proof}

It follows that:

\begin{theorem}
\label{divclAt}
The divisor class group of $A_t$ is free of rank $1$,
$$
\Cl(A_t)\iso \ZZ
$$
with basis $\cl(\qq)$.
Furthermore we have $\cl(\pp)=-t \cl(\qq)$.
\end{theorem}

\begin{proof} By \cite[Thm.\ 7.1]{F} and \ref{subintAt} there is an exact
sequence:
$$
0\to \ZZ\cl(\pp) \to \Cl(A_t) \to \Cl(V_t)\to 0.
$$
It is well-known that $\Cl(V_t)$ is isomorphic to $\ZZ/t\ZZ$ and
that it is spanned by the class of the prime ideal $P$ generated
by the elements of the  form $x_{11}\mu T$ where $\mu$ is a
monomial in the $x_{ij}$ of degree $t-1$. Note that $\cl(\pp)$ is
a torsion free element in $\Cl(A_t)$. This is because
$\pp^{(j)}=\{ F \in A_t : \gamma_2(F)\geq j\}$ cannot be
principal; in fact, it contains all the elements of the form $g^j$
where $g=f^tT^{t+1}$ and $f$ is a $(t+1)$-minor. Now fix a
$(t+1)$-minor $f$  and set $P_1=(f)S[T]\sect V_t$, that is  $P_1$ is
the ideal generated by all the elements of the form $f\mu T^2$
where $\mu$ is a monomial of degree $t-1$. Evidently $P_1$ is
isomorphic to $P$ and hence $\cl(P_1)$ generates $\Cl(V_t)$. But
$\qq$ is $P_1\sect A_t$ and hence $\cl(\qq)$ is the preimage of
$\cl(P_1)$ with respect to the map $\Cl(A_t) \to \Cl(V_t)$. It
follows that $\cl(\pp)$ and $\cl(\qq)$ generate $\Cl(A_t)$.  By
evaluating the function $\gamma_2$ one shows that
$$\qq^{(t)}\sect \pp=(g)A_t.$$ Hence $\cl(\pp)=-t\cl(\qq)$.
Consequently $\Cl(A_t)$
is generated by $\cl(\qq)$, and $\Cl(A_t)$ is torsion free.
\end{proof}

\section{Canonical modules of $\ini(\Rees)$ and $\ini(A_t)$} \label{canoini}

The goal of this section is to describe the canonical modules of
the semigroup rings  $\ini(\Rees)$ and $\ini(A_t)$. We know that
$\ini(\Rees)$ and $\ini(A_t)$ are normal (see \cite{BC})  and
hence their canonical modules are the vector spaces spanned by all
monomials represented by integral points in the relative interiors
of the corresponding cones.

To simplify notation we identify monomials of $S[T]$ with integral
points of $\RR^{mn+1}$. To a subset $G$ of the lattice
$[1,\dots,m]\times [1,\dots,n]$  we associate the ideal
$P_G=(x_{ij} : (i,j) \not\in G)$ of $S$ and a linear form $\ell_G$
on $\RR^{mn}$ defined by $\ell_G(x_{ij})=1$ if $(i,j) \not\in G$ and
$\ell_G(x_{ij})=0$ otherwise.

The initial ideal $\ini(I_i)$ of $I_i$ is the square free monomial
ideal generated by the initial terms (i.e.\ main diagonal
products) of the $i$-minors. Therefore $\ini(I_i)$ is the
Stanley-Reisner ideal of a simplicial complex $\Delta_i$. Let
${\bf F}_i$ be the set of the facets of $\Delta_i$. Then
$\ini(I_i)=\Sect_{F \in {\bf F}_i} P_F $. The elements of ${\bf
F}_i$  are described by Herzog and Trung \cite{HT} in terms of families of
non-intersecting paths. It turns out that $\Delta_i$ is a pure
(even shellable) simplicial complex. We have shown in \cite{BC1}
that
$$
\ini(I_i^{(k)})=\Sect_{F \in {\bf F}_i} P_F^k. \eqno{(2)}
$$
and that
$$
\ini(I_t^k)=\Sect_{i=1}^t \Sect_{F\in  {\bf F}_i} P_F^{k(t+1-i)}.
\eqno{(3)}
$$

For every $F\in  {\bf F}_i$  we extend the linear form $\ell_F$ to
the linear form $L_F$ on  $\RR^{mn+1}$ by setting
$L_F(T)=-(t+1-i)$. Then the equations (2) and (3) imply

\begin{lemma}  \label{iRmon}
A monomial $N$ belongs to $\ini(\Rees)$ if and only if it has
non-negative exponents and $L_F(N)\geq 0$ for every $F \in  {\bf
F}_i$ and $i=1,\dots,t$.
\end{lemma}

This is the description of the semigroup in terms of linear
homogeneous inequalities. It follows from the general theory
(see Bruns and Herzog \cite[Ch.~6]{BH}) that the canonical module
$\omega(\ini(\Rees))$
of $\ini(\Rees)$ is the semigroup ideal of $\Rees$ generated by
the monomials $N$ with all exponents $\geq 1$ and $L_F(N)\geq 1$
for every $F \in  {\bf F}_i$ and $i=1,\dots,t$.

Let $\XX$ denote the product of all the variables $x_{ij}$ with
$(i,j)\in [1,\dots,m]\times [1,\dots,n]$. The canonical module
$\omega(\ini(\Rees))$ has a description in terms of $\XX$ and the
functions $\rho_i$:

\begin{lemma} \label{caninRe}
The canonical module  $\omega(\ini(\Rees))$ of $\ini(\Rees)$ is
the ideal
$$
\{ F \in S[T] : \XX T \mid F \text{ in } S[T]   \text{ and }
\rho_i(F)\geq 1  \text{ for every } i=1,\dots,t\}
$$
of $\ini(\Rees)$.
\end{lemma}

\begin{proof}
Let $N=MT^k$ be a monomial, where $M$ is a monomial in the
variables $x_{ij}$. Then, for a given $i$, $L_F(N)\geq 1$ for
every $F \in {\bf F}_i$ if and only if  $\ell_F(M)\geq k(t+1-i)+1$
for every $F \in  {\bf F}_i$.  By (2) this is equivalent to $M\in
\ini(I_i^{(k(t+1-i)+1)})$ which in turn is equivalent to
$\rho_i(M)\geq k(t+1-i)+1$. Summing up,  $L_F(N)\geq 1$ for every
$F \in  {\bf F}_i$ and $i=1,\dots,t$ if and only if $\rho_i(N)\geq
1$ for every $i=1,\dots,t$.
\end{proof}

Similarly the canonical module $\omega(\ini(A_t))$ has a
description in terms of the function $\rho_2$:

\begin{lemma}
\label{caninAt} The canonical module $\omega(\ini(A_t))$ of
$\ini(A_t)$ is the ideal
$$
\{ F \in V_t : \XX T \mid F \text{ in } S[T] \text{ and }
\rho_2(F)\geq 1 \}
$$
of $\ini(A_t)$.
\end{lemma}

For later application we record the following lemma. Its part (2)
asserts that $\XX$ is a ``linear'' element for the functions
$\rho_i$.

\begin{lemma}
\label{rhoXX} \leavevmode\kern5pt
\begin{itemize}
\item[(1)] $\rho_i(\XX)=(m-i+1)(n-i+1)$.
\item[(2)] Let $M$ be any monomial in the $x_{ij}$'s. Then
$\rho_i(\XX M)=\rho_i(\XX)+\rho_i(M)$ for every $i=1,\dots,\min(m,n)$.
  \end{itemize}
\end{lemma}

\begin{proof} Let $M$ be a monomial in the $x_{ij}$'s.
We know that $\rho_i(M)\geq k$ if and only if  $M\in
\ini(I_i^{(k)})$. From  equation $(2)$ we may deduce that
$$
\rho_i(M)=\inf\{ \ell_F(M) :  F\in {\bf F}_i\}.
$$
Note that  $\Delta_i$ is a pure simplicial complex of dimension equal to
the dimension of the determinantal ring defined by $I_i$ minus $1$. It
follows that $\ell_F(\XX)=(m-i+1)(n-i+1)$ for every facet $F$ of $\Delta_i$.
In particular, $\rho_i(\XX)=(m-i+1)(n-i+1)$.

Since  $\ell_F(NM)\ge \ell_F(N)+\ell_F(M)$ for all monomials $N,M$ and for
every $F$, we have $\rho_i(MN)\geq \rho_i(N)+\rho_i(M)$.
Conversely, let $G$ be a facet of $\Delta_i$ such that
$\rho_i(M)=\ell_G(M)$. Then $\ell_G(\XX M)= \ell_G(\XX)+
\ell_G(M)=\rho_i(\XX)+\rho_i(M)$. Hence $\rho_i(MN)\leq
\rho_i(N)+\rho_i(M)$, too.
\end{proof}

\section{The deformation lemma}
We have been able to identify the canonical modules of $\ini(\Rees)$ and
$\ini(A_t)$ and we would like to use this information to describe the
canonical modules of $\Rees$ and $A_t$. To this end we need a deformation
lemma:

\begin{lemma}
\label{deform1} Let $R=K[x_1,\dots,x_n]$ be a polynomial ring
equipped with a term order $\tau$ and with a grading induced by
positive weights $\deg(x_i)=v_i$. Let $B$ be a finitely generated
$K$-subalgebra of $R$ generated by homogeneous polynomials and $J$
be a homogeneous ideal of $B$. Denote by $\ini(B)$ and $\ini(J)$
the initial algebra and the initial ideal of $B$ and $J$
respectively. Then we have:
\begin{itemize}
\item[(1)] If $\ini(B)$ is finitely generated and $\ini(B)/\ini(J)$ is
Cohen-Macaulay, then $B/J$ is Cohen-Macaulay.
\item[(2)] If $\ini(B)$ is finitely generated and Cohen-Macaulay and
$\ini(J)$ is the canonical module of $\ini(B)$ (up to shift) then
$B$ is Cohen-Macaulay and $J$ is the canonical module of $B$ (up
to the same shift).
\end{itemize}
\end{lemma}

\begin{proof}
(1) Let $f_1,\dots,f_k$ be a Sagbi basis of $B$ and
$g_1,\dots,g_h$ a Gr\"obner basis of $J$. We may assume that these
polynomials are monic and homogeneous. Consider the presentation
$K[y_1,\dots,y_k]/I\iso B$ of $B$ obtained by mapping $y_i$ to
$f_i$ and the  presentation $K[y_1,\dots,y_k]/I_1\iso \ini(B)$ of
$\ini(B)$ obtained by mapping $y_i$ to $\ini(f_i)$. Let
$h_1,\dots,h_p$ a system of binomial generators for the toric
ideal $I_1$, say $h_i=y^{a_i}-y^{b_i}$. Then for each $i$ we have
expressions $ f^{a_i}-f^{b_i}=\sum \lambda_{ij}f^{c_{ij}}$ with
$\lambda_{ij}\in K\setminus\{0\}$ and
$\ini(f^{c_{ij}})<\ini(f^{a_i})$ for every $i,j$. It is known that
the polynomials
$$
y^{a_i}-y^{b_i}-\sum \lambda_{ij}y^{c_{ij}} \eqno{(4)}
$$
generate $I$. For each $g_i$ we may take a presentation
$g_i=f^{d_i}+\sum \delta_{ij}f^{d_{ij}}$ with  $\delta_{ij}\in
K\setminus\{0\}$ and $\ini(f^{d_{ij}})<\ini(f^{d_i})$ for all
$i,j$. Then, by construction, the preimage of $J$ in $K[y]/I$  is
generated by the elements
$$
y^{d_i}+\sum \delta_{ij}y^{d_{ij}},\eqno{(5)}
$$
and the preimage of $\ini(J)$ in $K[y]/I_1$ is generated by the
$y^{d_i}$. Hence $B/J$ is the quotient of $K[y]$ defined by the
ideal $H$ that is generated by the polynomials (4) and (5), and
$\ini(B)/\ini(J)$ is the quotient of $K[y]$ defined by the ideal
$H_1$ that is generated by the polynomials $h_i$ and $y^{d_i}$.

If we can find  a positive weight $w$ on $K[y]$ such that
$\ini_w(H)=H_1$, then there is $1$-parameter flat family with
special fiber $\ini(B)/\ini(J)$ and general fiber $B/J$. This
implies that $B/J$ is Cohen-Macaulay, provided $\ini(B)/\ini(J)$
is. Let us define $w$. First consider a positive weight $\alpha$
on $K[x]$ such that $\ini_\alpha(f^{c_{ij}})<\ini_\alpha(f^{a_i})$
for every $i,j$ and $\ini_\alpha(f^{d_{ij}})<\ini_\alpha(f^{d_i})$
for every $i$ and $j$. That such an $\alpha$ exists is a
well-known property of monomial orders; for instance, see
Sturmfels \cite[Proof of Cor.\ 1.11]{S}. Then we define $w$ as the
``preimage" of $\alpha$ in the sense that we put
$w(y_i)=\alpha(\ini(f_i))$. It is clear, by construction, that the
initial forms of the polynomials (4) and (5) with respect to $w$
are exactly the $h_i$ and the $y^{d_i}$.

This proves that $\ini_w(H)$ contains $H_1$. But $H_1$ and $H$
have the same Hilbert function by construction, and $H$ and
$\ini_w(H)$ have the same Hilbert function because they have the same
initial ideal   if we refine $w$ to a term order; for
instance see \cite[Prop.\ 1.8]{S}. Here we consider Hilbert
functions with respect to the original graded structure induced by
the weights $v_i$. It follows that $\ini_w(H)=H_1$ and we are
done.

(2) That $B$ is Cohen-Macaulay follows from \cite{CHV}. Since $B$
is a Cohen-Macaulay positively graded $K$-algebra which is a
domain, to prove that $J$ is the canonical module of $B$ it
suffices to show that $J$ is a maximal Cohen-Macaulay module whose
Hilbert series satisfies the relation
$H_J(t)=(-1)^dt^kH_B(t^{-1})$ for some integer $k$ where $d=\dim
B$ \cite[Thm.\ 4.4.5,Cor.\ 4.4.6]{BH}.

The relation $H_J(t)=(-1)^dt^kH_B(t^{-1})$ holds since by
assumption the corresponding relation holds for the initial
objects and Hilbert series do not change by taking initial terms.
So it is enough to show that $J$ is a maximal Cohen-Macaulay
module. But $\ini(J)$ is a height $1$ ideal since it is the
canonical module \cite[Prop.\ 3.3.18]{BH}, and hence also $J$  has
height $1$. Therefore it suffices to show that $B/J$ is a
Cohen-Macaulay ring. But this follows from $(1)$ since
$\ini(B)/\ini(J)$ is Cohen-Macaulay (it is even Gorenstein)
\cite[Prop.\ 3.3.18]{BH}.
\end{proof}

\section{The canonical classes of $\Rees$ and $A_t$}

In this section we will describe the canonical  modules  of $\Rees$ and of
$A_t$ and determine the canonical classes.

Assume for  simplicity that $m\leq n$.  Let us consider a product
of minors $D$ such that $\ini(D)=\XX$ and
$\gamma_i(D)=\rho_i(\XX)$. Since  we have already computed
$\rho_i(\XX)$ (see \ref{rhoXX}), we can determine the shape of
$D$, which turns out to be $1^2,2^2,\dots,(m-1)^2, m^{(n-m+1)}$.
In other words, $D$ must be the product of $2$ minors of size $1$,
$2$ minors of size $2$, $\dots$, $2$ minors of size $m-1$ and
$n-m+1$ minors of size $m$. It is then not difficult to show that
$D$ is uniquely determined, the $1$-minors are $[m|1]$ and
$[1|n]$, the $2$-minors are $[m-1,m|1,2]$ and $[1,2|n-1,n]$ and so
on.

We have:

\begin{theorem}
\label{canmodRA}\leavevmode\kern5pt
\begin{itemize}
\item[(1)] The canonical module of $\Rees$ is the ideal
$$
J=\{ F\in S[T] :  DT \mid F \mbox{  in }  S[T] \mbox{ and }
\gamma_i(F)\geq 1 \mbox{ for  }  i=1,\dots,t\}.
$$
\item[(2)] The canonical module of $A_t$ is the ideal
$$
J_1=\{ F\in V_t  :  DT \mid F \mbox{  in }  S[T] \mbox{ and }
\gamma_2(F)\geq 1\}.
$$
\end{itemize}
\end{theorem}

\begin{proof} By virtue of \ref{deform1} it suffices to show that $\ini(J)$
and $\ini(J_1)$ are the canonical modules of $\ini(\Rees)$ and
$\ini(A_t)$ respectively. A description of the canonical modules
of $\ini(\Rees)$ and $\ini(A_t)$ has been given in Section
\ref{canoini}. Therefore it is enough to check that  $\ini(J)$ is
exactly the ideal described in \ref{caninRe} and $\ini(J_1)$ is
the ideal described in \ref{caninAt}. Note that we may write
$$
J=DT\{F\in S[T]  :  \gamma_i(F)\geq 1-\gamma_i(DT) \mbox{ for  }
i=1,\dots,t\}
$$
and
$$
J_1=DT\{ F\in S[T] :   \gamma_2(F)\geq 1-\gamma_2(DT)\} \sect V_t .
$$
Furthermore, by virtue of \ref{rhoXX},
$$
\omega(\ini(\Rees))=\XX T \{ F\in S[T]  :  \rho_i(F)\geq 1-\rho_i(\XX T)
\mbox{  for  }  i=1,\dots,t\}
$$
and
$$
\omega(\ini(A_t))=\XX T \bigl\{ F\in S[T]  :  \rho_2(F)\geq
1-\rho_2(\XX T)\bigr\} \sect V_t.
$$
Since, by the very definition of $D$, we have $\ini(D)=\XX$ and
$\gamma_i(D)=\rho_i(\XX)$, it suffices to show that
$$
\ini(\{ F \in S[T] :  \gamma_i(F)\geq j \})=\{F \in S[T] :
\rho_i(F)\geq j \}.
$$
But this has (essentially) been proved in \cite{BC}.
\end{proof}

Now we determine the canonical class of $\Rees$.

\begin{theorem}
\label{canclR}
The canonical class of $\Rees$ is given by
\begin{align*}
\cl(\omega(\Rees))&=\sum_{i=1}^t \Bigl(2-(m-i+1)(n-i+1)+t-i\Bigr)
\cl(P_i)\\
&=\cl(I_t\Rees)+\sum_{i=1}^t (1-\height I_i)\cl(P_i)
\end{align*}
\end{theorem}

\begin{proof} The second formula for $\omega(\Rees)$ follows
from the first, since $\height I_i=(m-i+1)(n-i+1)$ and
$\cl(I_t\Rees)=\sum_{i=1}^t (t-i+1)\cl(P_i)$ by \ref{primdecRees}.

We have seen that
$$
\omega(\Rees)= DT\bigl\{F\in S[T]  :  \gamma_i(F)\geq
1-\gamma_i(DT) \mbox{ for  } i=1,\dots,t\bigr\}.
$$
We can get rid of  $DT$ and obtain a representation of
$\omega(\Rees)$ as a fractional ideal, namely,
$$
\omega(\Rees)= \bigl\{F\in S[T] : \gamma_i(F)\geq 1-\gamma_i(DT)
\mbox{ for  } i=1,\dots,t\bigr\}.
$$
It follows that $\omega(\Rees)=\sect P_i^{(1-\gamma_i(DT))}$.  As
$\gamma_i(DT)=(m-i+1)(n-i+1)-(t+1-i)$, we are done.
\end{proof}

For $A_t$ the situation is slightly more difficult since $D\notin
A_t$ in general. Therefore we need an auxiliary lemma:

\begin{lemma} \label{clqj}
Let $f_j$ be a $j$-minor of $X$ and $\qq_j=(f_j)S[T]\sect A_t$.
Then $\qq_j$ is a height $1$ prime ideal of $A_t$ and
$\cl(\qq_j)=(j-t)\cl(\qq)$.
\end{lemma}

\begin{proof}
Note that $\qq_{t+1}=\qq$ by definition.  Let $\Delta$ be a
product of minors. We say that $\Delta$ has \emph{tight  shape} if
its degree is divisible by $t$ and it has exactly $\deg(\Delta)/t$
factors. In other words, $\Delta$ has tight shape if
$\gamma_1(\Delta T^k)=0$ and $\gamma_2(\Delta T^k)=0$ where
$k=\deg(\Delta)/t$. Let $\Delta$ be a product of minors with tight
shape and $k=\deg(\Delta)/t$. Note that
$$
(\Delta)S[T]\sect  A_t=(\Delta T^k)A_t.
$$
Now fix $j\leq t$ and set $k=t-j+1$; the product
$\Delta=f_jf_{t+1}^{t-j}$ has  tight shape and hence $\qq_j\sect
\qq_{t+1}^{(t-j)}=(\Delta T^k)A_t$. Note that,  for obvious
reasons, $\qq_j$ does not contain $\qq$. It follows that $\qq_j$
is a prime ideal of  height $1$  and that  $\cl(\qq_j)=(j-t)\cl(\qq)$.  Now
take $j>t$
and set $k=j-t+1$; the product $\Delta=f_jf_{t-1}^{j-t}$ has tight
shape and hence, as above, we conclude that $\qq_j$ is prime of
height $1$ and $\cl(\qq_j)=(t-j)\cl(q_{t-1})$. Since we know
already that $\cl(q_{t-1})=-\cl(\qq)$, we are done.
\end{proof}

Now we can prove

\begin{theorem}
\label{canclAt}
The canonical class of $A_t$ is given by
$$
\cl(\omega(A_t))=(mn-tm-tn)\cl(\qq).
$$
\end{theorem}

\begin{proof}  Assume that $m\leq n$. Set $W=(DT)S[T]\sect   A_t$.  We have
seen that
$$
\omega(A_t)= W \sect \pp.
$$
Consequently
$$
\cl(\omega(A_t))= \cl(W)+\cl(\pp)=\cl(W)-t\cl(\qq).
$$
Note that $W$ can be written as the intersection of $(D)S[T]\sect
A_t$  and $(T)S[T]\sect A_t$. But the latter is the irrelevant
maximal ideal of $A_t$,  whence $W=(D)S[T]\sect A_t$. Further
$(D)S[T]\sect A_t$ can be written as an intersection of ideals
$\qq_j$. Taking into consideration the shape of $D$ and
\ref{clqj}, we have
$$
\cl(W)=\biggl(\sum_{j=1}^{m-1} 2(j-t)+
(n-m+1)(m-t)\biggr)\cl(\qq).
$$
Summing up, we get the desired result.
\end{proof}

As a corollary we have
\begin{theorem}
\label{GorAt}
The ring $A_t$ is Gorenstein if and only if  $mn=t(m+n)$.
\end{theorem}

Note that we assume $1<t<\min(m,n)$ and $m\neq n$ if
$t=\min(m,n)-1$ in the theorem.   We have observed in the
introduction that $A_t$ is indeed Gorenstein (and even factorial)
in the remaining cases.

\begin{remark}
It is also possible to derive Theorem \ref{canclAt} from Theorem
\ref{canclR} by a suitable generalization of \cite[(8.10)]{BV}.
\end{remark}

\section{Algebras of minors of generic  Hankel matrices}
\label{Hankel}

Let $S$ be the  polynomial ring $K[x_1,\dots,x_n]$ over some field
$K$. Choose a Hankel matrix $X$ with distinct entries
$x_1,\dots,x_n$; this means that  $X$ is an $a\times b$ matrix
$(x_{ij})$  with $x_{ij}=x_{i+j-1}$ and $a+b-1=n$. Let $I_t$ be
the ideal generated by the minors of size $t$ of $X$. It is known
that $I_t$ does not depend on the size of the matrix $X$
(provided, of course,  $X$ contains $t$-minors); it depends only
on $t$ and $n$. For a given $n$ it follows that $t$ may vary from
$1$ to $m$, where $m=[(n+1)/2]$  is  the integer part of
$(n+1)/2$. All the properties of the generic determinantal ideals
that we have used in the previous sections to study their algebras
of minors hold also for the determinantal ideals of Hankel
matrices. This has been shown in Conca \cite{C}. In particular:
\smallskip

\noindent (1)  the symbolic powers $I_t^{(k)}$ and the primary
decomposition of $I_t^k$ are described in terms of the
$\gamma$-functions;
\smallskip

\noindent (2)  the initial algebras are described in terms of the
corresponding functions for monomials, the $\rho$-functions;
\smallskip

\noindent (3) the $\rho$-functions can be described in terms of
the linear forms associated to the facets of the simplicial
complex $\Delta_i$  of $\ini(I_i)$; moreover $\Delta_i$ is a pure
simplicial  \smallskip complex;

\noindent (4) the determinantal ideal $I_t$ is prime and $\height
I_t=n-2t+2$; in particular
the minors of $X$ are irreducible polynomials,
\smallskip

\noindent (5) the product of all the variables is a linear element
for the $\rho$ functions, that is, if $M$ is a monomial and $\XX$
is the product of all the variables then
$\rho_i(M\XX)=\rho_i(M)+\rho_i(\XX)$.
\smallskip

Now fix a number $t$, $1\leq t\leq m=[(n+1)/2]$. Denote by $\Rees$
the Rees algebra of $I_t$ and by $A_t$ the algebra generated by
the minors of size $t$. These algebras have been studied in
\cite[Section 4]{C}. It turns out  that $\Rees$ and $A_t$ are
Cohen-Macaulay normal domains and that the dimension of $A_t$ is
$n$, unless $t=m$. In the latter case, the minors are
algebraically independent. The arguments of the previous sections
apply also to the present situation. One deduces the following
theorem:

\begin{theorem}
\label{Han1}\leavevmode\kern5pt
\begin{itemize}
\item[(1)] If $t<m$, then the divisor class group $\Cl(\Rees)$ of $\Rees$ is
free of rank $t$ with basis  $\cl(P_1),\dots, \cl(P_t)$.  Here $P_i$ is the
prime ideal of $\Rees$ associated with the valuation $\gamma_i$ and
$P_i \sect S=I_i$. The
canonical class is
\begin{align*}
\cl(\omega(\Rees))&=\sum_{i=1}^t \Bigl(-n+t+i\Bigr)
\cl(P_i)\\
&=\cl(I_t\Rees)+\sum_{i=1}^t (1-\height I_i)\cl(P_i)
\end{align*}
\item[(2)] If $t=m$ and $n$ is even, then $\Cl(\Rees)$ is free of rank
$1$ with basis $\cl(P)$ where $P=I_m\Rees$. In this case $\Rees$
is a complete intersection and  $\cl(\omega(\Rees))=0$.
\item[(3)] If $t=m$ and $n$ is odd, then $I_m$ is principal.
\end{itemize}
\end{theorem}

To prove the theorem one argues like in the generic case. The only
difference is that the product of minors $D$ such that
$\ini(D)=\XX$ and $\gamma_i(D)=\rho_i(\XX)$ now has shape $m,m-1$
if $n$ is odd, and shape $m,m$ if $n$ is even.

For the algebra $A_t$ we have:

\begin{theorem}
\label{Han2}\leavevmode\kern5pt
\begin{itemize}
\item[(1)] Assume that  $1<t<m$ and $n$ is even if $t=m-1$.  Then the
divisor class group $\Cl(A_t)$ of $A_t$ is  free of rank $1$ with
basis $\cl(\qq)$. Here $\qq$ is the prime ideal of $A_t$  defined
as $\qq=(f)S[T]\sect A_t$ where $f$ is a minor of size $t+1$. The
canonical class is $\cl(\omega(A_t))=(n-3t)\cl(\qq)$.
\item[(2)] If $t=m-1$ and $n$ is odd, then $A_t$ is isomorphic to the
coordinate ring $\Grass(m-1,m+1)$ of the Grassmann variety of the
subspaces of dimension $m-1$ in a vector space of dimension $m+1$.
In particular $A_t$ is factorial.
\item[(3)] If $t=m$ or $t=1$, then $A_{t}$ is a polynomial ring.
\end{itemize}
\end{theorem}

The only assertion which still needs an argument is $(2)$. To this
end note that $I_t$ is the ideal of maximal minors of the Hankel
matrix of size $t\times (n+1)-t$. This induces a surjective map
from $\Grass(t, n+1-t)$ to $A_t$. If we take $n$ odd and $t=m-1$,
then the dimension of $\Grass(t, n+1-t)$ is $n$, which is also the
dimension of $A_t$. Hence the map $\Grass(t, n+1-t)\to A_t$ is an
isomorphism.

\begin{theorem}
\label{Han3} The ring $A_t$ is Gorenstein if and only if one of
the following conditions holds:
\begin{itemize}
\item[(1)] $3t=n$,
\item[(2)] $t=m-1$ and $n$ is odd,
\item[(3)] $t=m$ or $t=1$.
\end{itemize}
\end{theorem}

Despite of the analogy, let us point out that there are some
important differences between the Rees algebras $\Rees$ and the
algebras $A_t$ for generic matrices and those for Hankel matrices.
In the Hankel case the algebras $\Rees$  and  $A_t$ are defined by
(Gr\"obner bases of) quadrics, one has $\ini(I_t^k)=\ini(I_t)^k$
for all $t$ and $k$, and all the results hold over an arbitrary
field, no matter what the characteristic is.


\begin{thebibliography}{HSV}

\bibitem[Bo]{Bo} N. Bourbaki, {\em  Alg\`ebre commutative, Chap. I--IX}.
Hermann, Masson, 1961--1983.

\bibitem[B]{B} W. Bruns, {\em Algebras defined by powers of
determinantal ideals}, J. of Algebra 142 (1991), 150--163.

\bibitem[BC1]{BC} W. Bruns and A. Conca, {\em KRS
and powers of determinantal ideals}, Compositio Math. 111 (1998),
111--122.

\bibitem[BC2]{BC1} W. Bruns and A. Conca, {\em KRS
and   determinantal ideals},  preprint 2000.

\bibitem[BH]{BH}  W. Bruns and  J. Herzog, {\em Cohen-Macaulay rings},
Cambridge Studies in Advanced Mathematics, 39,  Cambridge University Press,
Cambridge, 1993.

\bibitem[BST]{BST}  W. Bruns, A. Simis, and N. V. Trung,
{\em  Blow-up of straightening-closed ideals in ordinal Hodge
algebras}, Trans. Amer. Math. Soc. 326 (1991), 507--528.

\bibitem[BV]{BV} W. Bruns and U. Vetter, {\em Determinantal rings},
Lect. Notes Math. 1327, Springer 1988.

\bibitem[C]{C} A.Conca, {\em Straightening law and powers of determinantal
ideals of Hankel matrices}, Adv. Math. 138 (1998), 263-292.

\bibitem[CHV]{CHV} A. Conca, J. Herzog, and G. Valla, {\em Sagbi bases
and application to blow-up algebras}, J. Reine Angew. Math. {474} (1996),
113--138.

\bibitem[DEP]{DEP} C. De Concini, D. Eisenbud, and C. Procesi, {\em Young
diagrams and determinantal varieties}, Invent. math. 56 (1980),
129--165.

\bibitem[DRS]{DRS} P. Doubilet, G.C. Rota, and J. Stein, {\em On the
foundations of combinatorial theory: IX, Combinatorials methods in
invariants theory},
Studies in Applied Mathematics LIII (1974), 185--216.

\bibitem[EH]{EH}  D. Eisenbud and C. Huneke, {\em Cohen-Macaulay Rees algebras
and their specialization}, J. Algebra 81 (1983), 202--224.

\bibitem[F]{F} R. Fossum,  {\em The divisor class group of a Krull domain},
Ergebnisse der Mathematik und ihrer Grenzgebiete, Band 74.
Springer-Verlag, New York-Heidelberg, 1973.

\bibitem[HSV1]{HSV} J. Herzog, A. Simis, and  W. Vasconcelos,   {\em On
the canonical module of the Rees algebra and the associated graded
ring of an ideal}, J. Algebra 105 (1987), 285--302.

\bibitem[HSV2]{HSV1} J. Herzog, A. Simis, and  W. Vasconcelos, {\em
Arithmetic of normal Rees algebras}, J. Algebra 143 (1991),
269--294.

\bibitem[HV]{HV} J. Herzog and W. Vasconcelos,
{\em On the divisor class group of Rees-algebras}, J. Algebra 93
(1985), 182--188.

\bibitem[HT]{HT} J. Herzog and N. V. Trung, {\em Gr\"obner bases and
multiplicity of determinantal and pfaffian ideals}, Adv. in Math.
96 (1992), 1--37.

\bibitem[ST]{ST} A. Simis and N. V. Trung,  {\em The divisor class group
of ordinary and symbolic blow-ups}, Math. Z. 198 (1988), 479--491.

\bibitem[St]{St}  R. P. Stanley, {\em Enumerative combinatorics, Vol. 2},
Cambridge University Press 1999 .

\bibitem[S]{S}  B. Sturmfels, {\em Gr\"obner bases and convex polytopes}, AMS
University Lecture Series, Vol. 8, Providence RI, 1995.

\bibitem[V]{V} W. V. Vasconcelos, {\em Arithmetic of blowup algebras},  London
Mathematical Society Lecture Note Series, 195,  Cambridge University Press,
Cambridge, 1994.


\end{thebibliography}
  \end{document}